\documentclass[aip,
amsmath,amssymb,
preprint,%
numerical,
]{revtex4-1}
\usepackage{graphicx}
\usepackage{dcolumn}
\usepackage{bm}
\usepackage{hyperref}
\usepackage[utf8]{inputenc}
\usepackage[T1]{fontenc}
\usepackage{mathptmx}
\usepackage{etoolbox}
\newcommand{\RNum}[1]{\uppercase\expandafter{\romannumeral #1\relax}}
\usepackage{array,xcolor,colortbl}
\usepackage{multirow}
\makeatletter
\def\@email#1#2{%
	\endgroup
	\patchcmd{\titleblock@produce}
	{\frontmatter@RRAPformat}
	{\frontmatter@RRAPformat{\produce@RRAP{*#1\href{mailto:#2}{#2}}}\frontmatter@RRAPformat}
	{}{}
}%
\makeatother

\begin{document}
	
	\title{Effect of oblique irradiation on the onset of thermal-phototactic-bioconvection in an isotropic scattering algal suspension}
	
	\author{S. K. Rajput}
	\altaffiliation[Corresponding author: E-mail: ]{shubh.iiitj@gmail.com.}
	\author{M. K. Panda}%
	\affiliation{$^1$ Department of Mathematics, PDPM Indian Institute of Information Technology Design and Manufacturing, Jabalpur 482005, India.
	}%
	

	\begin{abstract}
		
	In this study, our focus is mainly to check the effect of light scattering on the onset of thermal-phototactic-bioconvection in an algal suspension where the suspension is illuminated by the collimated oblique irradiation from above while simultaneously applying heating or cooling from below. We conduct a numerical investigation into the linear stability of a suspension containing phototactic algae, focusing particularly on how the angle of incidence of oblique collimated irradiation influences the system. Our solutions reveal a transition of the most unstable mode from a stationary to an overstable state, or vice versa, under certain parameter configurations as the angle of incidence varies. Additionally, we frequently observe oscillatory instabilities in cases where the upper surface is rigid, particularly as the angle of incidence increases within the suspension.
		
	\end{abstract}
	
	
	\maketitle
	
	
	\section{INTRODUCTION}
	
	Bioconvection is an intriguing occurrence marked by fluid movement driven by mobile microorganisms within the fluid, a phenomenon extensively examined in previous research~\cite{20platt1961,21pedley1992,22hill2005,23bees2020,24javadi2020}. Typically found in water, these microorganisms display a distinctive behavior of swimming upwards due to their higher density relative to the surrounding medium. The emergence of specific patterns in bioconvection is intricately tied to the swimming behavior of these microorganisms. However, pattern formation isn't solely dictated by their upward movement or greater density; it's also influenced by their response to various stimuli known as "taxes," which include gravitaxis, chemotaxis, phototaxis, gyrotaxis, and thermotaxis. Phototaxis involves movement in response to light, while thermotaxis entails directed movement in response to temperature changes, observed in both organisms and cells. The complex interplay of these environmental cues leads to the diverse array of patterns observed in bioconvection. The term "thermal photactic bioconvection" specifically denotes a subset of bioconvection where light and temperature gradients are crucial in governing the motion and behaviors of microorganisms, particularly motile algae. This phenomenon explores the unique interaction between phototaxis and thermotaxis in these microorganisms. Exploring thermal photactic bioconvection offers researchers an opportunity to investigate the combined effects of these environmental cues on the movement and distribution of aquatic microorganisms.	
	
	The examination of bioconvection has been extensively pursued across various domains, particularly focusing on the interplay between thermal and phototactic factors in suspensions of microorganisms. Kuznetsov~\cite{51kuznetsov2005thermo} delved into bio-thermal convection within a suspension of oxytactic microorganisms, while Alloui et al. \cite{52alloui2006stability}  investigated suspensions of mobile gravitactic microorganisms. Nield and Kuznetsov \cite{53nield2006onset} employed linear stability analysis to probe the onset of bio-thermal convection in gyrotactic microorganism suspensions, and Alloui et al.~\cite{54alloui2007numerical} scrutinized the influence of bottom heating on the onset of gravitactic bioconvection in a square enclosure. Taheri and Bilgen \cite{55taheri2008thermo} explored the effects of bottom heating or cooling in a vertically oriented cylinder with stress-free sidewalls. Kuznetsov \cite{56kuznetsov2011bio} developed a theoretical framework for bio-thermal convection in suspensions containing two species of microorganisms. Saini et al. \cite{57saini2018analysis} investigated bio-thermal convection in suspensions of gravitactic microorganisms, while Zhao et al.~\cite{57zhao2018linear} utilized linear stability analysis to examine bioconvection stability in suspensions of randomly swimming gyrotactic microorganisms heated from below.
	
    \begin{figure}[!htbp]
		\centering
		\includegraphics[width=14cm]{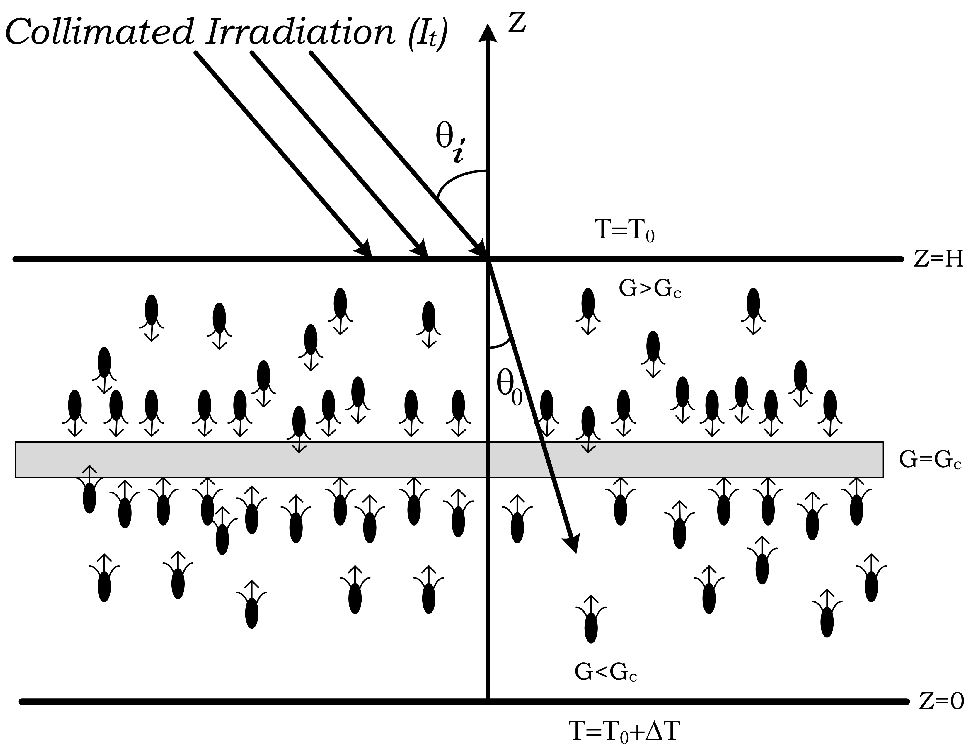}
		\caption{\footnotesize{Formation of the sublayer within the suspension due to oblique collimated irradiation.}}
		\label{fig1}
	\end{figure}
	
	In the realm of phototactic bioconvection, Vincent and Hill~\cite{12vincent1996} laid the foundation with groundbreaking research, examining the effects of collimated irradiation on an absorbing (non-scattering) medium. Ghorai and Hill~\cite{10ghorai2005} extended this inquiry into the behavior of phototactic algal suspensions in two dimensions, without considering scattering effects. Ghorai et al.~\cite{7ghorai2010} and Ghorai and Panda~\cite{13ghorai2013} investigated light scattering, both isotropic and anisotropic, under normal collimated irradiation. Panda and Ghorai~\cite{14panda2013} proposed a model for an isotropically scattering medium in two dimensions, yielding results divergent from those of Ghorai and Hill~\cite{10ghorai2005} due to the inclusion of scattering effects. Panda and Singh~\cite{11panda2016} explored phototactic bioconvection in two dimensions, confining a non-scattering suspension between rigid sidewalls. Panda et al.~\cite{15panda2016} examined the impact of diffuse irradiation, combined with collimated irradiation, in an isotropic scattering medium, while Panda~\cite{8panda2020} investigated an anisotropic medium. Considering natural environmental conditions where sunlight strikes the Earth's surface at oblique angles, Panda et al.~\cite{16panda2022} studied the effects of oblique collimated irradiation on the onset of phototactic bioconvection. In a recent investigation, Panda and Rajput~\cite{41rajput2023} explored the impacts of diffuse irradiation along with oblique collimated irradiation on a uniformly scattering suspension. Rajput and Panda~\cite{rajput2024effect} investigated the effect of scattered/diffuse flux on the onset of phototactic bioconvection in the absence of collimated flux. 
	
   The extensive investigation into thermal bioconvection has been of paramount importance due to its biological significance. However, the current body of literature lacks exploration into thermal phototactic bioconvection within an isotropic scattering medium which is illuminated by oblique collimated flux. This research presents a fresh perspective on this phenomenon by elucidating the significant influence of angle of incidence with thermal effects on the phototactic behavior of organisms within such a scattering medium. Integrating thermal effects into phototactic bioconvection not only increases the complexity of this natural phenomenon but also unveils new dimensions in comprehending the intricate interplay between temperature, light, and biology.

   The structure of this article is organized as follows: Firstly, the mathematical formulation of the problem is presented, followed by the derivation of a fundamental (equilibrium) solution. Subsequently, a small disturbance is introduced to the equilibrium system, and the linear stability problem is obtained through the application of linear perturbation theory, followed by numerical solution methods. The model's results are then presented, and finally, the implications and findings of the model are thoroughly discussed. This systematic approach facilitates a comprehensive exploration of the interaction between thermal gradients, light, and biological responses within the context of isotropic scattering media.
	

\section{GEOMETRY OF THE PROBLEM}
	
	Consider the motion in a dilute suspension of phototactic algae within a layer of finite depth $H$. Here, an oblique collimated irradiation
	illuminates the suspension from above and strikes it at a fixed off-normal angle $\theta_i$: We choose a rectangular Cartesian coordinate system where the yz-plane is the plane of incidence for the oblique collimated irradiation (see Fig.~\ref{fig1}). The angle of refraction $\theta_0$ in which the collimated beam propagates across the water is determined by using the Snell’s law. For simplicity, the effects of scattering by algae have been neglected in the present study. Also we consider the temperature differences is gentle enough to not pose any lethal threat to the microorganisms. Furthermore, the phototactic behavior, encompassing cell orientation and swimming speed, remains unaffected by these thermal variations. At the bottom boundary $(z=0)$, the temperature is kept at $T_0+\Delta T$, while at the top boundary $(z=H)$, it remains uniformly at $T_0$. These temperature conditions persist uniformly across the suspension.
	
	\begin{figure}[!htbp]
		\centering
		\includegraphics[width=14cm]{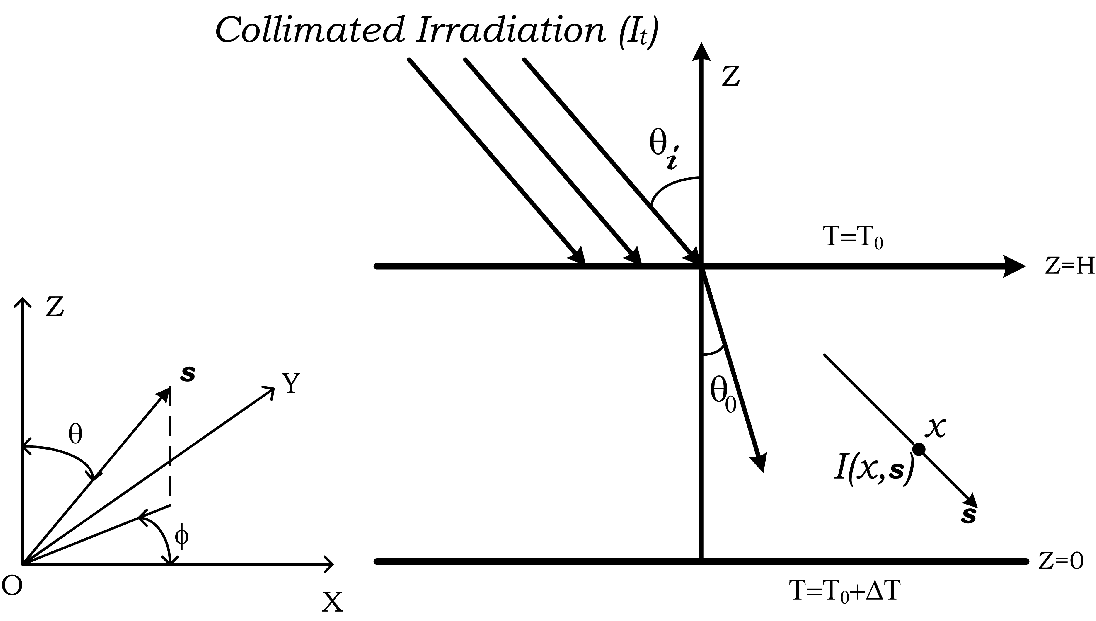}
		\caption{\footnotesize{Oblique collimated irradiation on the upper surface of non-scattering algal suspension.}}
		\label{fig2}
	\end{figure}

\section{PHOTOTAXIS IN A SCATTERING SUSPENSION}
	
	Let $I(\boldsymbol{x},\boldsymbol{s})$ represents the radiation intensity propagating in the (unit) direction $\boldsymbol{s}$ at a position $\boldsymbol{x}$ across the algal suspension, where $\boldsymbol{x}$ is measured relative to a rectangular cartesian coordinate system with the $z$ axis vertically up (see Fig.~\ref{fig2}). We assume here that the medium across the algal suspension to be non-scattering similar to Vincent and Hill. To calculate light intensity profiles, the radiative transfer equation (hereafter referred to as RTE) is given by
	
	\begin{equation}\label{1}
	\frac{dI(\boldsymbol{x},\boldsymbol{s})}{ds}+(a+\sigma_s)I(\boldsymbol{x},\boldsymbol{s})=\frac{\sigma_s}{4\pi}\int_{0}^{4\pi}I(\boldsymbol{x},\boldsymbol{s'})p(\boldsymbol{s},\boldsymbol{s'})d\Omega'.
	\end{equation}
	
	Here, $\Omega$ is the solid angle and $p(\boldsymbol{s},\boldsymbol{s'})$ is the scattering phase function. For simplicity, the proposed work assumes that the scattering by algae across the suspension is isotropic (i.e., $p=1$).
	
	At location $\boldsymbol{x}_l=(x,y,H)$ the light intensity is given by
	
	\begin{equation*}
	I(\boldsymbol{x}_l,\boldsymbol{s}) = I_t\delta(\boldsymbol{s},\boldsymbol{s}_0).
	\end{equation*}
	
	where $I_t$ is the magnitude of the oblique collimated flux in the direction $\boldsymbol{s}_0=\sin(\pi-\theta_0)\cos\phi\hat{x}+\sin(\pi-\theta_0)\sin\phi\hat{y}+\cos(\pi-\theta_0)\hat{z}$.

  We assume $a=\alpha n(x)$ and $\sigma_s=\beta n(x)$, and the RTE now
  becomes 

\begin{equation}\label{2}
\frac{dI(\boldsymbol{x},\boldsymbol{s})}{ds}+(\alpha+\beta)nI(\boldsymbol{x},\boldsymbol{s})=\frac{\beta n}{4\pi}\int_{0}^{4\pi}I(\boldsymbol{x},\boldsymbol{s'})d\Omega'.
\end{equation}

The total intensity $G(\boldsymbol{x})$ is defined as

\begin{equation*}
G(\boldsymbol{x})=\int_0^{4\pi}I(\boldsymbol{x},\boldsymbol{s})d\Omega,
\end{equation*}

and the light flux $\boldsymbol{q}(\boldsymbol{x})$ is

\begin{equation}\label{3}
\boldsymbol{q}(\boldsymbol{x})=\int_0^{4\pi}I(\boldsymbol{x},\boldsymbol{s})\boldsymbol{s}d\Omega.
\end{equation}

The average swimming velocity is given by Hill and Hader

\begin{equation*}
\boldsymbol{W}_c=W_c<\boldsymbol{p}>.
\end{equation*}

Let $<\boldsymbol{p}>$ denotes the average swimming direction, 15 and it is given
by
\begin{equation}\label{4}
<\boldsymbol{p}>=-M(G)\frac{\boldsymbol{q}}{|\boldsymbol{q}|},
\end{equation}

Here, where $T(G)$ is a photoresponse curve such that

\begin{equation*}
T(G)=\left\{\begin{array}{ll}\geq 0, & \mbox{if }~~ G\leq G_{c},\\< 0, & \mbox{if }~~G>G_{c}.  \end{array}\right. 
\end{equation*}

\section{THE GOVERNING BIOCONVECTIVE SYSTEM}

Reminiscent to the previous continuum models on bioconvection, the volume fraction of the algae is small as the suspension is dilute and the interaction between the cells is negligible. Let the volume and density of each cell be denoted by $\vartheta$ and $\rho + \Delta\rho$, respectively, where $\Delta\rho$ is the difference in density between cell and water and $0<\Delta\rho\ll\rho$: Let the cell concentration and average velocity of all the material in a small volume $\delta\vartheta$ be defined as $n$ and $\boldsymbol{u}$, respectively. Then, suspension incompressibility implies

\begin{equation}\label{5}
	\boldsymbol{\nabla}\cdot \boldsymbol{u}=0.
\end{equation}

The Navier–Stokes equations under the Boussinesq approximation become
\begin{equation}\label{6}
	\rho\left(\frac{D\boldsymbol{u}}{D t}\right)=-\boldsymbol{\nabla} P_e+\mu{\nabla}^2\boldsymbol{u}-nv g\Delta\rho\hat{\boldsymbol{z}}-\rho g(1-\beta(T-T_0))\hat{\boldsymbol{z}}.
\end{equation}

Here, $D/Dt$ denotes the convective derivative, $P_e$ is the excess hydrostatic pressure, and $\mu$ is the suspension (or water) viscosity.

The equation governing the conservation of algae is given by
\begin{equation}\label{7}
	\frac{\partial n}{\partial t}=-\boldsymbol{\nabla}\cdot \boldsymbol{F},
\end{equation}

where $\boldsymbol{F}$ is the cell flux and it is written as $\boldsymbol{F}=n(\boldsymbol{u}+W_c<\boldsymbol{p}>)-\boldsymbol{D}\cdot\boldsymbol{\nabla} n$.

The thermal energy equation is expressed as
\begin{equation}\label{8}
	\rho c\big[\frac{\partial T}{\partial t}+\boldsymbol{\nabla}\cdot(\boldsymbol{u}T)\big] =\alpha{\boldsymbol{\nabla}}^2 T,
\end{equation}

where, $\rho c$ is the volumetric heat capacity of water, and $\alpha$ is the thermal conductivity of water.

The appropriate boundary conditions are

\begin{subequations}
	\begin{equation}\label{9a}
		\boldsymbol{u}\cdot\hat{\boldsymbol{z}}=\boldsymbol{F}\cdot\hat{\boldsymbol{z}}=0\qquad \text{at }~~~~ z=0,H.
	\end{equation}
	The rigid boundary conditions at $z=0,H$ are
	\begin{equation}\label{9b}
		\boldsymbol{u}\times\hat{\boldsymbol{z}}=0\qquad \text{at }~~~~ z=0,H.
	\end{equation}
	while the free boundary conditions at $z=0,H$ are
	\begin{equation}\label{9c}
		\frac{\partial^2}{\partial z^2}(\boldsymbol{u}\cdot\hat{\boldsymbol{z}})=0\qquad \text{at }~~~~ z=0,H.
	\end{equation}
	For temperature
	\begin{equation}\label{9d}
		T=T_0+\Delta T\qquad \text{at } ~~~~z=0,
	\end{equation}
	\begin{equation}\label{9e}
		T=T_0\qquad \text{at }~~~~ z=H.
	\end{equation}
\end{subequations}

The boundary intensities are given as
\begin{equation*}
	I(x,y,z=H,\theta,\phi)=I_t\delta(\boldsymbol{s},\boldsymbol{s}_0)~~~~~~~~(\pi/2\leq\theta\leq\pi).
\end{equation*}
\begin{equation*}
	I(x,y,z=0,\theta,\phi)=0~~~~~~~~(0\leq\theta\leq\pi/2).
\end{equation*}

To derive dimensionless governing equations, we normalize lengths with the layer depth $H$, scale time using the diffusive time scale $H^2/\alpha_f$, and express bulk fluid velocity in terms of $\alpha_f/H$. Pressure is scaled by $\mu \alpha_f/H^2$, cell concentration is normalized with $\bar{n}$ (the mean concentration), and temperature is scaled by $(T-T_0)/\Delta T$. The resulting non-dimensional bioconvection equations are as follows

\begin{equation}\label{10}
	\boldsymbol{\nabla}\cdot\boldsymbol{u}=0,
\end{equation}
\begin{equation}\label{11}
	P_r^{-1}\left(\frac{D\boldsymbol{u}}{D t}\right)=-\nabla P_{e}+\nabla^{2}\boldsymbol{u}-R_bn\hat{\boldsymbol{z}}-R_m\hat{\boldsymbol{z}}+R_TT\hat{\boldsymbol{z}},
\end{equation}
\begin{equation}\label{12}
	\frac{\partial{n}}{\partial{t}}=-\boldsymbol{\nabla}\cdot[\boldsymbol{n{\boldsymbol{u}}+\frac{1}{Le}nV_{c}<{\boldsymbol{p}}>-\frac{1}{Le}{\boldsymbol{\nabla}}n}],
\end{equation}
and
\begin{equation}\label{13}
	\frac{\partial T}{\partial t}+\boldsymbol{\nabla}\cdot(\boldsymbol{u}T) =\boldsymbol{\nabla}^2 T.
\end{equation}

The boundary conditions after scaling become

\begin{subequations}
	\begin{equation}\label{14a}
		\boldsymbol{u}\cdot\hat{\boldsymbol{z}}=\big[\boldsymbol{n{\boldsymbol{u}}+\frac{1}{Le}nV_{c}<{\boldsymbol{p}}>-\frac{1}{Le}{\boldsymbol{\nabla}}n}\big]\cdot\hat{\boldsymbol{z}}=0\qquad \text{at } z=0,1.
	\end{equation}
	The rigid boundary conditions at $z=0,1$ become
	\begin{equation}\label{14b}
		\boldsymbol{u}\times\hat{\boldsymbol{z}}=0\qquad \text{at } z=0,1.
	\end{equation}
	while free boundary conditions at $z=1$ become
	\begin{equation}\label{14c}
		\frac{\partial^2}{\partial z^2}(\boldsymbol{u}\cdot\hat{\boldsymbol{z}})=0\qquad \text{at }~~~~ z=1.
	\end{equation}	
	For temperature
	\begin{equation}\label{14d}
		T=1\qquad \text{at }~~~~ z=0,
	\end{equation}
	\begin{equation}\label{14e}
		T=0\qquad \text{at }~~~~ z=1.
	\end{equation}
\end{subequations}

The RTE [see Eq. (1)] becomes

\begin{equation}\label{15}
	\frac{dI(\boldsymbol{x},\boldsymbol{s})}{ds}+\kappa nI(\boldsymbol{x},\boldsymbol{s})=\frac{\sigma n}{4\pi}\int_{0}^{4\pi}I(\boldsymbol{x},\boldsymbol{s'})d\Omega',
\end{equation}

where $\kappa=(\alpha+\beta)\Bar{n}H$ and $\sigma=\beta\Bar{n}H$ are the extinction and scattering coefﬁcients after scaling. The single scattering albedo is deﬁned as $\omega=\sigma/\kappa$. Equation (23) in terms of $\omega$ can be written as

\begin{equation}\label{16}
	\frac{dI(\boldsymbol{x},\boldsymbol{s})}{ds}+\kappa nI(\boldsymbol{x},\boldsymbol{s})=\frac{\omega\kappa n}{4\pi}\int_{0}^{4\pi}I(\boldsymbol{x},\boldsymbol{s'})d\Omega'.
\end{equation}

Here, $\omega\in[0,1]$ and $\omega=0 (\omega=1)$ represents a purely absorbing (scattering) medium. The top and bottom intensities after scaling
become

In the form of direction cosine, RTE becomes:

\begin{equation}\label{17}
	\nu_1\frac{dI}{dx}+\nu_2\frac{dI}{dy}+\nu_3\frac{dI}{dz}+\kappa nI(\boldsymbol{x},\boldsymbol{s})=\frac{\omega\kappa n}{4\pi}\int_{0}^{4\pi}I(\boldsymbol{x},\boldsymbol{s'})d\Omega',
\end{equation}

where $\xi,\eta$ and $\nu$ are the direction cosines in x, y and z direction. In dimensionless form, the intensity at boundaries becomes:

\begin{equation*}
	I(x,y,z=1,\theta,\phi)=I_t\delta(\boldsymbol{s},\boldsymbol{s}_0)~~~~~~~~(\pi/2\leq\theta\leq\pi).
\end{equation*}
\begin{equation*}
	I(x,y,z=0,\theta,\phi)=0~~~~~~~~(0\leq\theta\leq\pi/2).
\end{equation*}

\section{THE EQUILIBRIUM STATE}

In the basic state, we set
\begin{equation}\label{18}
\boldsymbol{u}=0,~n=n_p(z),~\text{and}~ I=I_p(z,\theta)
\end{equation}
in Eqs.~(\ref{10})–(\ref{13}) and (\ref{17}) and solve.

In the basic state, the total intensity $G_p$ and light ﬂux $\boldsymbol{q}_p$ become

\begin{equation*}
	G_p=\int_0^{4\pi}I_p(z,\theta)d\Omega,\quad
	\boldsymbol{q}_p=\int_0^{4\pi}I_p(z,\theta)\boldsymbol{s}d\Omega.
\end{equation*}

The governing equation for $I_p$ becomes

\begin{equation}\label{19}
\frac{dI_p}{dz}+\frac{\kappa n_pI_p}{\nu}=\frac{\omega\kappa n_p}{4\pi\nu}G_p(z),
\end{equation}
where $\nu_1=\sin\theta\cos\phi$, $\nu_2=\sin\theta\sin\phi$, and $\nu_3=\cos\theta$ are the direction cosines in $x$, $y$, and $z$ directions, respectively.

Now, $I_p$ can be divided as $I_p=I_p^c+I_p^d$, where $I_p^c$ is the collimated part after attenuation and $I_p^d$ is the diffused part. The collimated part, $I_p^c$, is governed by

\begin{equation}\label{20}
\frac{dI_p^c}{dz}+\frac{\kappa n_pI_p^c}{\nu}=0,
\end{equation}

subject to the boundary condition,

\begin{equation*}
	I_p^c( z=1, \theta) =I_t\delta(\boldsymbol{s}-\boldsymbol{ s}_0),~~~~~~ (\pi/2\leq\theta\leq\pi),
\end{equation*}

and $I_p^d$ is governed by

\begin{equation}\label{21}
\frac{dI_p^d}{dz}+\frac{\kappa n_pI_p^d}{\nu_3}=\frac{\omega\kappa n_p}{4\pi\nu_3}G_p(z),
\end{equation}

where boundary conditions are

\begin{equation*}
	I_p^d( z=1, \theta) =0,~~~~~~ (\pi/2\leq\theta\leq\pi),
\end{equation*}

Now, the total intensity is decomposed as the sum of a collimated and diffused one, i.e.,

\begin{equation}\label{22}
G_p=G_p^c+G_p^d,
\end{equation}	
where	
\begin{equation}\label{23}
G_p^c=\int_0^{4\pi}I_p^c(z,\theta)d\Omega=I_t\exp\left(\frac{\kappa\int_1^z n_p(z')dz'}{\cos\theta_0}\right),
\end{equation}

\begin{equation}\label{24}
G_p^d=\int_0^{4\pi}I_p^d(z,\theta)d\Omega.
\end{equation}
The Lambert–Beer law, i.e., $G_p=G_p^c$, is obtained if we neglect scattering.	

In terms of $\tau=\kappa\int_z^1 n_p(z')dz'$ as the optical thickness, the total intensity after scaling, i.e., $\Upsilon(\tau)=G_p(\tau)/I_t$ becomes

\begin{equation}\label{25}
\Upsilon(\tau) = e^{-\tau/\cos\theta_0}+\frac{\omega}{2}\int_0^\kappa \Upsilon(\tau')E_1(|\tau-\tau'|)d\tau',
\end{equation}
which is a Fredholm integral equation, and $E_n(x)$ is the exponential integral of order $n$. Equation (\ref{25}) is solved using the method of substraction of singularity	

The light ﬂux in the basic state becomes		\begin{equation*}
	\boldsymbol{q_p}=\int_0^{4\pi}\left(I_p^c+I_p^d\right)\boldsymbol{s}d\Omega=-I_t\cos\theta_0\exp\left(\frac{\kappa\int_1^z n_p(z')dz'}{\cos\theta_0}\right)\hat{\boldsymbol{z}}+\int_0^{4\pi}I_p^d(z,\theta)\boldsymbol{s}d\Omega.
\end{equation*}

Now, $\boldsymbol{q}_s=-q_p\hat{\boldsymbol{z}}$, where $q_p=|\boldsymbol{q_p}|$ as $I_p^d$ is independent of
$\phi$. Now, the average swimming orientation is given by

\begin{equation*}
	<\boldsymbol{p_s}>=M_p\hat{\boldsymbol{z}},
\end{equation*}

where $M_p=M(G_p)$.

The concentration $n_p(z)$ satisfies

\begin{equation}\label{26}
\frac{dn_p}{dz}=V_cM_pn_p,
\end{equation}

associated with

\begin{equation}\label{27}
\int_0^1n_p(z)dz=1.
\end{equation}

Here, $M_p = M(G)$ at $G = G_p$. The steady light intensity $G_p$ at a height $z$ $(0 \leq z \leq 1)$ is defined as $T_p(z)$, satisfying

\begin{equation}\label{28}
\frac{d^2T_p}{dz^2}=0.
\end{equation}

The boundary conditions (\ref{14d}) and (\ref{14e}) lead to

\begin{equation}\label{29}
T_p(z)=1-z.
\end{equation}

The boundary value problem via Eqs. (\ref{25})-(\ref{29}) can be solved by a shooting method.

In our study, we consider Chlamydomonas to estimate the parameter ranges. The range of the angle of incidence $\theta_i$ for the proposed model is such that $0\leq\theta_i\leq 80$. The optical depth $\kappa\in[0.25:1]$ for a 0.5 cm deep suspension as given in the study by Ghorai and Panda. Hence, the calculated scaled swimming speed are $V_c=10$ and $V_c=20$,respectively. The magnitude of the oblique collimated ﬂux ($I_t$) at the top surface has been assumed to be equal to unity.
	
	
\section{Linear stability of the problem}
The basic state, i.e., Eq.~(\ref{18}) is perturbed via a small perturbation of amplitude $\epsilon$ (where $0 < \epsilon \ll 1$) such as
\begin{equation}\label{30}
\begin{pmatrix}
\boldsymbol{u}\\n\\T\\<\boldsymbol{I}>
\end{pmatrix}
=
\begin{pmatrix}
0\\n_p\\T_p\\<\boldsymbol{I}_p>
\end{pmatrix}
+\epsilon
\begin{pmatrix}
\boldsymbol{u}_1\\n_1\\T_1\\<\boldsymbol{I}_1>
\end{pmatrix}
+O(\epsilon^2),
\end{equation}

where $\boldsymbol{u}_1=(u_1,v_1,w_1)$. The perturbed variables are substituted into Eqs.~(\ref{10})–(\ref{13}). Next, we linearize about the basic state and collect the $O(\epsilon)$ terms. So that

\begin{equation}\label{31}
\boldsymbol{\nabla}\cdot \boldsymbol{u}_1=0,
\end{equation}

\begin{equation}\label{32}
P_r^{-1}\left(\frac{\partial \boldsymbol{u_1}}{\partial t}\right)=-\boldsymbol{\nabla} P_{e}+\nabla^{2}\boldsymbol{u}_1-R_bn_1\hat{\boldsymbol{z}}+R_TT_1\hat{\boldsymbol{z}},
\end{equation}

\begin{equation}\label{33}
\frac{\partial{n_1}}{\partial{t}}+\frac{1}{Le}V_c\boldsymbol{\nabla}\cdot(<\boldsymbol{p_s}>n_1+<\boldsymbol{p_1}>n_p)+w_1\frac{dn_p}{dz}=\frac{1}{Le}\boldsymbol{\nabla}^2n_1,
\end{equation}

\begin{equation}\label{34}
\frac{\partial{T_1}}{\partial t}-w_1\frac{dT_S}{dz}=\boldsymbol{\nabla}^2T_1.
\end{equation}

If $G=G_p+\epsilon G_1+O(\epsilon^2)$,
then the steady collimated total intensity is  perturbed as\\ $I_t\exp\left(\frac{\kappa\int_1^z (n_p(z')+\epsilon n_1+O(\epsilon^2))dz'}{\cos\theta_0}\right)$
and after simpliﬁcation, we get

\begin{equation}\label{35}
G_1^c=I_t\exp\left(\frac{\kappa\int_1^z n_p(z')dz'}{\cos\theta_0}\right)\left(\frac{\kappa\int_1^z n_1 dz'}{\cos\theta_0}\right),
\end{equation}

Similarly, $G_1^d$ is given by

\begin{equation}\label{36}
G_1^d=\int_0^{4\pi}I_1^d(\boldsymbol{ x},\boldsymbol{ s})d\Omega,
\end{equation}

Similarly,

\begin{equation}\label{37}
\boldsymbol{q}_1^c=-I_t\exp\left(\frac{\kappa\int_1^z n_p(z')dz'}{\cos\theta_0}\right)\left(\frac{\kappa\int_1^z n_1 dz'}{\cos\theta_0}\right)\cos\theta_0\hat{z}
\end{equation}

\begin{equation}\label{38}
q_1^d=\int_0^{4\pi}I_1^d(\boldsymbol{x},\boldsymbol{s})\boldsymbol{s}d\Omega.
\end{equation}

By using a similar process, we find the perturbed swimming direction

\begin{equation}\label{39}
<\boldsymbol{p_1}>=G_1\frac{dM_p}{dG}\hat{\boldsymbol{z}}-M_p\frac{\boldsymbol{q_1}^H}{\boldsymbol{q_p}},
\end{equation}

Note that the last term becomes zero for non-scattering algal suspension illuminated by oblique collimated ﬂux.

Equation~(\ref{39}) is inserted into Eq.~(\ref{33}) and after simpliﬁcation,
we get

\begin{equation}\label{40}
\frac{\partial{n_1}}{\partial{t}}+\frac{1}{L_e}V_c\frac{\partial}{\partial z}\left(T_pn_1+n_pG_1\frac{dT_p}{dG}\right)-\frac{1}{L_e}V_cn_p\frac{T_p}{q_p}\left(\frac{\partial q_1^x}{\partial x}+\frac{\partial q_1^y}{\partial y}\right)+w_1\frac{dn_p}{dz}=\frac{1}{L_e}\nabla^2n_1.
\end{equation}

We eliminate $u_1$, $v_1$ of $\boldsymbol{u_1}$ and $P_e$ from the perturbed governing bioconvection system, so that Eqs.~(\ref{32}), (\ref{33}), (\ref{34}) and (\ref{40}) are reduced into two equations for $w_1$ and $n_1$. Thus

\begin{equation}\label{41}
\begin{pmatrix}
w_1\\n_1\\T_1
\end{pmatrix}
=
\begin{pmatrix}
W(z)\\\Theta(z)\\T(z)
\end{pmatrix}
+\exp{[\sigma t+i(k_1x+k_2y)]},
\end{equation} 

$W(z)$, $\Theta(z)$, and $T(z)$ represent the variations in the $z$ direction, while $k_1$ and $k_2$ are the horizontal wavenumbers.

In terms of direction cosines $(\nu_1, \nu_2, \nu_3)$, $I_1^d$ can be written as
\begin{equation}\label{42}
\nu_1\frac{\partial I_1^d}{\partial x}+\nu_2\frac{\partial I_1^d}{\partial y}+\nu_3\frac{\partial I_1^d}{\partial z}+\kappa n_pI_1^d=\frac{\omega\kappa}{4\pi}(n_pG_1^c+n_pG_1^d+G_pn_1)-\kappa n_1I_p,
\end{equation}
subject to the boundary conditions

\begin{subequations}
	\begin{equation}\label{43a}
	I_1^d(x, y, z=1, \nu_1, \nu_2, \nu_3) =0,~~~where~~~ (\pi/2\leq\theta\leq\pi,~~0\leq\phi\leq 2\pi), 
	\end{equation}
	\begin{equation}\label{43b}
	I_1^d(x, y, z=0,\nu_1, \nu_2, \nu_3) =0,~~~where~~~ (0\leq\theta\leq\pi/2,~~0\leq\phi\leq 2\pi). 
	\end{equation}
\end{subequations}

$I_1^d$ can be expressed as

\begin{equation*}
	I_1^d=\Psi_1^d(z,\nu_1,\nu_2,\nu_3)\exp{(\sigma t+i(k_1x+k_2y))}. 
\end{equation*}

From Eqs.~(\ref{35}) and (\ref{36}), we get

\begin{equation}\label{44}
G_1^c=\left[I_t\exp\left(-\int_z^1 \kappa n_p(z')dz'\right)\left(\int_1^z\kappa n_1 dz'\right)\right]\exp{(\sigma t+i(k_1x+k_2y))}=\mathcal{G}_1^c(z)\exp{(\sigma t+i(k_1x+k_2y))},
\end{equation}
and 
\begin{equation}\label{45}
G_1^d=\mathcal{G}_1^d(z)\exp{(\sigma t+i(k_1x+k_2y))}= \left(\int_0^{4\pi}\Psi_1^d(z,\nu_1,\nu_2,\nu_3)d\Omega\right)\exp{(\sigma t+i(k_1x+k_2y))},
\end{equation}

where $\mathcal{G}_1(z)=\mathcal{G}_1^c(z)+\mathcal{G}_1^d(z)$.

Now $\Psi_1^d$ satisfies
\begin{equation}\label{46}
\frac{d\Psi_1^d}{dz}+\frac{(i(k_1\nu_1+k_2\nu_2)+\kappa n_p)}{\nu_3}\Psi_1^d=\frac{\omega\kappa}{4\pi\nu_3}(n_p\mathcal{G}_1+G_p\Theta)-\frac{\kappa}{\nu_3}I_p\Theta,
\end{equation}

with subject to

\begin{subequations}
	\begin{equation}\label{47a}
	\Psi_1^d( z=1, \nu_1, \nu_2, \nu_3) =0,~~~\text{where}~~~ (\pi/2\leq\theta\leq\pi,~~0\leq\phi\leq 2\pi), 
	\end{equation}
	\begin{equation}\label{47b}
	\Psi_1^d( z=0,\nu_1, \nu_2, \nu_3) =0,~~~\text{where}~~~ (0\leq\theta\leq\pi/2,~~0\leq\phi\leq 2\pi). 
	\end{equation}
\end{subequations}

Equation (\ref{46}) (i.e., an integrodifferential equation) is solved via a suitable iterative method.

Similarly from Eqs.~(\ref{37}) and (\ref{38}), we have

\begin{equation*}
	q_1^H=[q_1^x,q_1^y]=[A(z),B(z)]\exp{[\sigma t+i(k_1x+k_2y)]},
\end{equation*}
where
\begin{equation*}
	A(z)=\int_0^{4\pi}\Psi_1^d(z,\nu_1,\nu_2,\nu_3)\nu_1 d\Omega,\quad B(z)=\int_0^{4\pi}\Psi_1^d(z,\nu_1,\nu_2,\nu_3)\nu_2 d\Omega.
\end{equation*}
The linear stability equations become
\begin{equation}\label{48}
\left(\sigma P_r^{-1}+k^2-\frac{d^2}{dz^2}\right)\left( \frac{d^2}{dz^2}-k^2\right)W=R_bk^2\Theta(z)-R_Tk^2T(z),
\end{equation}
\begin{equation}\label{49}
\left(Le\sigma+k^2-\frac{d^2}{dz^2}\right)\Theta(z)+V_c\frac{d}{dz}\left(T_p\Theta+n_p\mathcal{G}_1\frac{dT_p}{dG}\right)-i\frac{V_cn_pT_p}{q_p}(k_1A+k_2B)=-Le\frac{dn_p}{dz}W,
\end{equation}

\begin{equation}\label{50}
\left(\frac{d^2}{dz^2}-k^2-\sigma\right)T(z)=\frac{dT_s}{dz}W(z),
\end{equation} 

with
\begin{equation}\label{51}
W=\frac{d^2W}{dz^2}=\frac{d\Theta}{dz}-V_cT_p\Theta-n_pV_c\mathcal{G}_1\frac{dT_p}{dG}=0,~~\text{at}~~~z=0,1,
\end{equation}
where $k=\sqrt{k_1^2+k_2^2}$.

Eq.~(\ref{49}) becomes (writing D = d/dz)
\begin{equation}\label{52}
D^2\Theta-\aleph_3(z)D\Theta-(Le\sigma+k^2\aleph_2(z))\Theta-\aleph_1(z)\int_1^z\Theta dz-\aleph_0(z)=LeDn_pW, 
\end{equation}
where
\begin{subequations}
	\begin{equation}\label{53a}
	\aleph_0(z)=V_cD\left(n_p\mathcal{G}_1^d\frac{dT_p}{dG}\right)-i\frac{V_cn_pT_p}{q_p}(k_1A+k_2B),
	\end{equation}
	\begin{equation}\label{53b}
	\aleph_1(z)=\kappa V_cD\left(n_pG_p^c\frac{dT_p}{dG}\right),
	\end{equation}
	\begin{equation}\label{53c}
	\aleph_2(z)=2\kappa V_c n_p G_p^c\frac{dT_p}{dG}+V_c\frac{dT_p}{dG}DG_p^d,
	\end{equation}
	\begin{equation}\label{53d}
	\aleph_3(z)=V_cT_p.
	\end{equation}
\end{subequations}
Now, consider
\begin{equation}\label{54}
\Phi(z)=\int_1^z\Theta(z')dz',
\end{equation}
Eq.~(\ref{48}),~(\ref{50}) and (\ref{52}) becoems
\begin{equation}\label{55}
D^4W-(2k^2+\sigma Le P_r^{-1})D^2W+k^2(k^2+\sigma Le P_r^{-1})W=-R_bk^2D\Phi+R_Tk^2T(z),
\end{equation}
\begin{equation}\label{56}
D^3\Phi-\aleph_3(z)D^2\Phi-(\sigma Le+k^2+\aleph_2(z))D\Phi-\aleph_1(z)\Phi-	\aleph_0(z)=Le Dn_pW, 
\end{equation}
\begin{equation}\label{57}
\left(D^2-k^2-\sigma\right)T(z)=DT_sW,
\end{equation}
with

\begin{equation}\label{58}
W=D^2W=D^2\Phi-\Gamma_2(z)D\Phi-V_cn_p\frac{dT_p}{dG}\mathcal{G}_1=0,~~~ \text{at} ~~~z=0,1,
\end{equation}
\begin{equation}\label{59}
T(z)=0,~~~at~~z=0,1,
\end{equation}
and
\begin{equation}\label{60}
\Phi(z)=0,~~~ \text{at}~~~z=1.
\end{equation}
	
	
\section{SOLUTION PROCEDURE}

Equations~(\ref{55}), (\ref{56}) and (\ref{57}) are solved using a scheme via
Newton–Raphson–Kantorovich (NRK) iterations~\cite{19cash1980} and the marginal stability curves are computed in the $(k,R)$-plane, where $R=(R_b, R_T)$. The most unstable solution is recognized as the pair $(k^c=2\pi/\lambda_c, R^c)$, where $R_c$ is the minimum value of $R$ on the neutral curve $R^{(n)}(k)$, with $n=1,2,3..$, and $\lambda_c$ is the pattern wavelength. The bioconvective solution is named as mode $n$ if it has $n$ convection cells arranged vertically one on another. The graph of points for which $Re(\sigma)=0$ is called a marginal stability curve (or neutral curve). Furthermore, if $Im(\sigma)\neq 0$ on such a curve, then the perturbation to basic state is called stationary (or oscillatory/overstable). In the case of oscillatory instability, a single oscillatory branch originates from the stationary branch at some $k\leq k_b$ (say).


\section{NUMERICAL RESULTS}

In the proposed work, the lower boundary is assumed to be rigid
and the upper boundary is stress-free as similar to bioconvection
experiments. We have systematically investigated the effect of angle of incidence $\theta_i$ by varying it between 0 and 80 (i.e., $0\leq\theta_i\leq80$),
keeping $P_r, I_t, G_c, V_c, \omega$, and $\kappa$ ﬁxed. We take here a discrete set of parameters in our study in contrast to the whole parameter domain. We ﬁx $P_r=5$, $I_t=1$; $V_c=10,$ 15, 20, $\kappa=0:5$, 1.0, and $\omega\in[0,1]$, respectively, throughout this section. Based on weak scattering and strong scattering by algae, the results obtained for a discrete parameter set are divided into two categories. Since we are investigating the effect of oblique collimated irradiation on the suspension with two variety of the top boundary (stress free and rigid). So we devide results into two sections.	
	

\subsection{WHEN TOP SURFACE IS STRESS FREE}
\begin{figure}[!htbp]
	\centering
	\includegraphics[width=16cm]{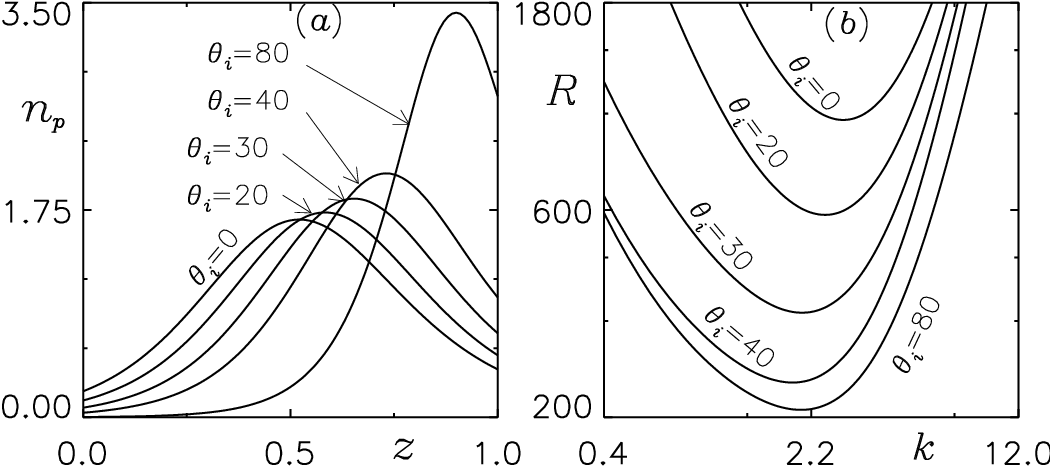}
	\caption{\footnotesize{(a) The profiles formed at basic state and (b) the corresponding marginal stability curves as $\theta_i$ is increased.}}
	\label{fig3}
\end{figure}
Fig.~\ref{fig3} illustrates the proﬁles
formed at the basic state and the corresponding marginal stability curves at the variation of $\theta_i$ , when the parameters $V_c=15$, $\kappa=0.5$; and $\omega=0.4$ are kept ﬁxed. We start with the case when $\theta_i=0$. In this case, the algal suspension is illuminated by a normal/vertical collimated solar radiation. To study the effects of collimated ﬂux on bioconvective instability, let the sublayer where the algae aggregate at basic state locates at around the mid-height of the suspension. In this instance, let the sublayer location be at some height $z=z_1$ inside the suspension, where the algae receive an optimal light via normal incidence, i.e., $G=G_c$. Let an off-normal collimated ﬂux be the illuminating source to the algal suspension and it is deﬁned through an angle of incidence, say $\theta_i=20$. In this case, let the sublayer location be at some height $z=z_2$ inside the suspension, where the algae receive an optimal light via off-normal incidence, i.e., $G=G_c$ : Then $z_1<z_2$ and this implies the width (thickness) of the upper stable layer overlying
the unstable layer decreases as compared to the case when $\theta_i=0$ (i.e., normal/vertical collimated irradiation). Thus, the buoyancy of the stable
layer, which tends to inhibit convective ﬂuid motions becomes less effective and the value of $R_c$ decreases than the case when $\theta_i=0$. When $\theta_i$ is increased further to 40, the sublayer of algae at the basic state is located around three-quarter height of the suspension. In this case, the non-oscillatory branch accommodates the most unstable mode of disturbance making the bioconvective solution to be stationary. Also,
both the critical Rayleigh number and wavenumber decrease due to the fact that the buoyancy of the upper stable zone becomes less effective against the bioconvective ﬂuid motions than the previous cases.
As $\theta_i$ increases further to 60, the base concentration proﬁle becomes steeper. This results in increment in the critical wavenumber in comparison to the previous case, but the critical Rayleigh number
decreases further. When $\theta_i$ reaches to 80, the steepness in the base concentration proﬁle increases and the thickness of the stable sublayer further decreases. Thus, both the critical Rayleigh number and
wavenumber increase (see Fig.~\ref{fig3}).

\subsection{WHEN TOP SURFACE IS RIGID}

	\begin{figure}[!htbp]
	\centering
	\includegraphics[width=16cm]{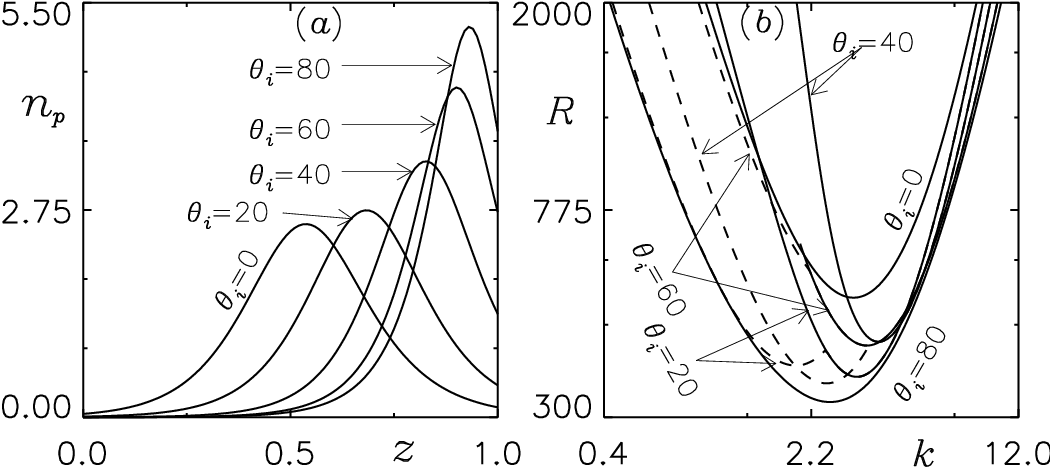}
	\caption{\footnotesize{(a) The profiles formed at basic state and (b) the corresponding marginal stability curves as $\theta_i$ is increased.}}
	\label{fig4}
    \end{figure}

In this section, we discuss the case when $V_c=15$, $\kappa=1$, and $\omega=0.4$. For a vertical/normal collimated solar irradiation, i.e., angle of incidence $\theta_i=0$, the sublayer of algae at basic state is located at around the mid-height of the chamber [see Fig.~\ref{fig4}(a)] and here the perturbation to the basic state becomes stationary.
When an off-normal/oblique collimated solar ﬂux is introduced as an illuminating source to the suspension via $\theta_i=20$; the sublayer of algae at the basic state is located at around $z=0.65$[see Fig.~\ref{fig4}(a)]. As explained in Sec. VIII, the width of the upper stable zone decreases as compared to the case of a vertical/normal collimated solar ﬂux, i.e., $\theta_i=0$. In this instance, the combination of buoyancy of the upper stable zone and positive phototaxis in the lower unstable zone can reinforce each other and inhibit the bioconvective ﬂow motions, whereas the negative phototaxis in the upper stable zone supports bioconvection. Eventually, the competition among these mechanisms resulted in
bifurcation of a single oscillatory branch from the stationary branch of the marginal stability curve at around $k\approx2$ and deﬁnes a graph for $k<2$ when $\theta_i=20$. However, the non-oscillatory
(stationary) branch accommodates the unstable mode of disturbance here and, thus, making the bioconvective solution to be stationary. At $\theta_i=40$, the sublayer of algae at basic state occurs at around $z= 0.8$ and a single oscillatory branch bifurcates from the stationary branch at around $k_b=3.4$ and the oscillatory branch retains the most unstable solution making the bioconvective solution to be overstable. Here, overstability occurs at $k_c=2.22$ and $R_c=329.53$. In this instance, the bifurcation is called Hopf (periodic) bifurcation. Stationary instability is observed as $\theta_i=60$ and 80 [see Fig.~\ref{fig4}].
	
	
	\section{Conclusion}
	
In this innovative model of thermal phototactic bioconvection driven by the combined effect of heat and oblique collimated irradiation, we delve into the occurrence of bio-thermal convection within a suspension comprising isotropic scattering phototactic algae. This suspension is illuminated by oblique collimated irradaitaion from top and heated from below. Our main aim is to examine the collective impacts of temperature and light scattering on the onset of thermo-phototactic bioconvection. We also investigate the linear instability analysis by utilizing the linear perturbation theory.

The sublayer where the cells accumulate at the base state shifts toward the top-half of the suspension at the increment in angle of incidence. Similarly, the value of maximum basic concentration increases at the increment in angle of incidence. Usually, the critical Rayleigh number at suspension instability decreases at the increment in angle of incidence for a weak-scattering algal suspension. At bioconvective instability, the perturbation to the basic steady state shifts from a stationary state to an oscillatory state or vice versa at the increment in angle of incidence/diffuse ﬂux for a weak-scattering suspension.

The initial pattern size (or critical wavelength) at suspension
instability increases with an increase (a decrease) in the angle of incidence (light intensity) for governing parameters. Thus, the experimental reports provided by Williams and Bees and the theoretical predictions via the proposed model are well matched. However, algae are predominantly photo-gyrotactic or photo-gravitactic or photo-gyro-gravitactic in nature to the best of our knowledge. Thus, we fail to include here a comparison between the proposed work with an appropriate quantitative study.

Finally, we note that algae in a natural environment absorb the
light via an illuminating source and scatter it in the forward direction across the suspension. Thus, a suitable modiﬁcation in the scattering phase function, i.e., $p(\boldsymbol{s},\boldsymbol{s'})$, may suggest another new and fascinating research direction related to phototactic bioconvection in a forward scattering algal suspension in the same vicinity.

	
	\nocite{*}
	\bibliography{Isotropic_oblique}
	
\end{document}